\documentclass{amsart}
\usepackage{amsfonts}
\usepackage{amsmath,amssymb}
\usepackage{amsthm}
\usepackage{amscd}
\usepackage{graphics}
\usepackage{graphicx}
\theoremstyle{definition}{
\newtheorem{Def}{{\rm Definition}}
\newtheorem{Ex}{{\rm Example}}
\newtheorem{Rem}{{\rm Remark}}

}
\theoremstyle{plain}
{

\newtheorem{Prop}{Proposition}
\newtheorem{Thm}{Theorem}
\newtheorem{MainThm}{Main Theorem}

}

\begin{document}
\title[Complex projective spaces have no special generic maps]{On the non-existence of special generic maps on complex projective spaces}
\author{Naoki Kitazawa}
\address{Institute of Mathematics for Industry, Kyushu University, 744 Motooka, Nishi-ku Fukuoka 819-0395, Japan\\
 TEL (Office): +81-92-802-4402 \\
 FAX (Office): +81-92-802-4405 \\
}
\email{n-kitazawa@imi.kyushu-u.ac.jp}
\urladdr{https://naokikitazawa.github.io/NaokiKitazawa.html}

\begin{abstract}
We prove the non-existence of {\it special generic} maps on complex projective space as our extended new result.
Simplest special generic maps are Morse functions with exactly two singular points on spheres, or Morse functions in Reeb's theorem, and canonical projections of unit spheres.
 
 Manifolds represented as connected sums of products of manifolds diffeomorphic to unit spheres admit such maps in considerable cases.

Real and complex projective spaces have been shown to admit no such maps in most cases by the author. This gives a complete answer for complex projective spaces as a corollary to a more general result, which is also our main result.


\end{abstract}


\maketitle

\section{Introduction.}
\label{sec:1}
  
Morse theory, or theory of Morse functions, gives fundamental and strong tools in singularity theory and its applications to geometry of manifolds.    
  
A {\it singular} point of a smooth map means a point in the manifold of the domain where the rank of the differential is smaller than both the manifolds of the domain and the target.  

\begin{Thm}
	A closed manifold admits a Morse function with exactly two singular points if and only if it is homeomorphic to a sphere in the case where the manifold is not $4$-dimensional. In the $4$-dimensional case, this holds if it is diffeomorphic to the unit sphere.
\end{Thm}

These functions are generalized as follows.

\begin{Def}
	A smooth map from an $m$-dimensional manifold with no boundary into an $n$-dimensional manifold with no boundanry is called a {\it special generic} map if $m \geq n \geq 1$ and for suitable coordinates, at each singular point it is represented by the form $(x_1,\cdots,x_m) \rightarrow (x_1,\cdots,x_{n-1},{\Sigma}_{j=1}^{m-n+1} {x_{n-1+j}}^2)$
\end{Def}

The definition seems to restrict the topology and the differntiable structure of the manifold strongly. This is true. 

Before explaining about this, we introduce fundamental terminologies and notation on smooth manifolds.
${\mathbb{R}}^k$ denotes the $k$-dimensional Euclidean space with the canonical differentiable structure. This is a simplest smooth manifold. This is also regarded as the Riemannian manifold endowed with the standard Euclidean metric.
For $x \in {\mathbb{R}}^k$, $||x||$ denotes the distance between $x$ and the origin $0 \in {\mathbb{R}}^k$ under the standard metric or the canonical Euclidean metric.
For a positive integer $k$, $D^k:=\{x \mid ||x|| \leq 1\} \subset {\mathbb{R}}^k$ denotes the $k$-dimensional unit disk, which is a $k$-dimensional compact, connected and smooth closed submanifold. $S^{k-1}:=\{x \mid ||x||=1\} \subset {\mathbb{R}}^{k}$ denotes the ($k-1$)-dimensional unit sphere, which is a ($k-1$)-dimensional closed smooth submanifold with no boundary. 
The $0$-dimensional unit sphere is a two-point set endowed with the discrete topology for $k=1$ ($k-1=0$). The ($k-1$)-dimensional unit sphere is connected for $k \geq 2$ ($k-1 \geq 1$).

Canonical projections of unit spheres are special generic. {\it Exotic} spheres, which are defined as manifolds which are homeomorphic to (unit) spheres and which are not diffeomorphic to unit spheres, do not admit such maps in considerable cases. \cite{calabi,saeki,saeki2, saekisakuma,saekisakuma2,wrazidlo} give related answers for example. 

On the other hands, for example, manifolds represented as connected sums of manifolds diffeomorphic to unit spheres admit these maps in considerable cases.
  
It is important and difficult problems to know the (non-)existence of special generic maps.  
  
For example, it seems to be natural to consider such problems for projective spaces, $2$, $3$ or $4$-dimensional closed manifolds and higher dimensional closed (and simply-connected) manifolds for example. In fact related results such as \cite{saeki,saeki2,saekisakuma,saekisakuma2,sakuma} have been obtained by Saeki and Sakuma since 1990. Recently Nishioka (\cite{nishioka}), Wrazidlo (\cite{wrazidlo,wrazidlo2,wrazidlo3}) and the author (\cite{kitazawa1,kitazawa2,kitazawa3,kitazawa4}) have obtained related results.
Our paper concerns such problems for complex projective spaces. The (non-)existence for real-projective spaces has been announced to be completely solved by the author. In short, existence is rare. 

\begin{MainThm}
	\label{mthm:1}
	Let $k \geq 2$ be an arbitrary integer.
	The $k$-dimensional complex projective space, which is a $2k$-dimensional closed and simply-connected smooth manifold, admits no special generic maps into ${\mathbb{R}}^{2k-1}$. Furthermore, by Theorem \ref{thm:4}, it admits no special generic maps into any Euclidean space.  
	\end{MainThm}

This was an open problem for an arbitrary odd integer $k \geq 3$ in \cite{kitazawa3}.
This follows from the following more general theorem. Terminologies and notation on homology groups and cohomology groups and rings are in section \ref{sec:3}. However, we suppose elementary knowledge on these notions.

Here we note that ${\mathbb{R}}^1$ is generally denoted by $\mathbb{R}$, that the ring of all integers (rational numbers) is denoted by $\mathbb{Z}$ (resp. $\mathbb{Q}$), that the relation $\mathbb{Z} \subset \mathbb{Q} \subset \mathbb{R}$ holds and that $\mathbb{Q}$ is the field of quotients (fractions) of the free commutative group $\mathbb{Z}$, whose rank is $1$. 
\begin{MainThm}
	\label{mthm:2}
	Let $k_1$ and $k_2$ be arbitrary integers greater than $1$.
	Let $M$ be a closed and simply-connected smooth manifold of dimension $m=k_1k_2$ satisfying the following conditions.
	\begin{enumerate}
		\item The rational cohomology group $H^{j}(M;\mathbb{Q})$ is isomorphic to $\mathbb{Q}$ if $0 \leq j \leq k_1k_2$ and $j$ is divisible by $k_1$ and it is the trivial group if $j$ is not as before. 
		\item There exists a generator $u_{0}$ of $H^{k_1}(M;\mathbb{Q})$ and the cup product ${\cup}_{j_1=1}^{j_2} u_0$, or the $j_2$-th power of $u_0$ in the rational cohomology ring $H^{\ast}(M;\mathbb{Q})$ is a generator of $H^{j_2k_1}(M;\mathbb{Q})$ for every integer $1 \leq j_2 \leq k_2-1$.
	\end{enumerate}
	Then $M$ admits no special generic maps into ${\mathbb{R}}^{m-1}$.  
\end{MainThm}

The second section is for fundamental properties of special generic maps. We also review results related to our present study more precisely.

The third section is devoted to expositions on Main Theorems \ref{mthm:1} and \ref{mthm:2}. In fact, a key ingredient is merely a slight generalization of a family of arguments used by the author in \cite{kitazawa2,kitazawa3,kitazawa4}. We close our paper by related examples and remarks. 

\section{Preliminaries.}
\label{sec:2}

Throughout the present paper, a smooth manifold homeomorphic to a unit sphere is called a {\it homotopy sphere}. A homotopy sphere which is (not) diffeomorphic to a unit sphere is a {\it standard} sphere (resp. an {\it exotic} sphere).

Every smooth manifold is well-known to admit the structure of a PL manifold canonically. In several scenes, we regard smooth manifolds as the canonical PL manifolds. More generally, topological manifolds are known to be regarded as CW complexes.

For a manifold or a polyhedron $X$, let $\dim X$ denote the dimension of $X$.

A {\it diffeomorphism} is a smooth map having no singular points being also a homeomorphism. A diffeomorphism from a given manifold $X$ onto the same manifold is said to be a {\it diffeomorphism} on $X$. The diffeomorphism group of the smooth manifold $X$ is the group of all diffeomorphisms on $X$. We also regard this as the topological space endowed with the so-called {\it Whitney $C^{\infty}$ topology}. {\it Whitney $C^{\infty}$ topologies} are natural topologies of sets of smooth maps between smooth manifolds. See \cite{golubitskyguillemin} for systematic expositions for example. 

We also need terminologies and notions on bundles. For systematic expositions, see \cite{milnorstasheff,steenrod}.
 
A bundle is said to be {\it smooth} if the fiber is a smooth manifold and the structure group is a subgroup of the diffeomorphism group. A {\it linear} bundle is a smooth bundle whose fiber is an Euclidean space, a unit sphere, or a unit disk and whose structure group consists of diffeomorphisms regarded as linear transformations naturally.

\begin{Prop}[E. g. \cite{saeki}]
\label{prop:1}
Let $m \geq n \geq 1$ be integers.
A special generic map $f:M \rightarrow N$ on an $m$-dimensional closed manifold $M$ into an $n$-dimensional manifold $N$ with no boundary always enjoys the following properties.
\begin{enumerate}
\item \label{prop:1.1}
There exist an $n$-dimensional compact manifold $W_f$, a smooth immersion $\bar{f}:W_f \rightarrow N$ and a smooth surjection $q_f:M \rightarrow W_f$ satisfying the relation $f=\bar{f} \circ q_f$.
\item \label{prop:1.2}
$q_f$ maps the singular set $S(f)$ of $f$ onto the boundary $\partial W_f \subset W_f$ as a diffeomorphism. Moreover, $S(f)$ is a smooth closed submanifold in $M$ with no boundary.
\item \label{prop:1.3}
 We have the following two bundles.
\begin{enumerate}
\item \label{prop:1.3.1}
For a suitable small collar neighborhood $N(\partial W_f) \subset W_f$, we have a linear bundle whose fiber is the {\rm(}$m-n+1${\rm )}-dimensional unit disk and whose projection is given by the composition of the map $q_f {\mid}_{{q_f}^{-1}(N(\partial W_f))}$ onto $N(\partial W_f)$ with the canonical projection to $\partial W_f$. 
\item \label{prop:1.3.2}
We have a smooth bundle over $W_f-{\rm Int}\ N(\partial W_f)$ whose fiber is an {\rm (}$m-n${\rm )}-dimensional standard sphere and whose projection is given by the restriction of $q_f$ to the preimage of $W_f-{\rm Int}\ N(\partial W_f)$.
\end{enumerate}
\end{enumerate}
\end{Prop}

Hereafter, $|X|$ denotes the size of a finite set $X$.

\begin{Ex}
\label{ex:1}
Let $m \geq n \geq 2$ be integers. Let $\{S^{k_j} \times S^{m-k_j}\}_{j \in J}$ be a family of finitely many products of two unit spheres where $k_j$ is an integer satisfying $1 \leq k_j \leq n-1$. We consider a connected sum of these $|J|$ manifolds in the smooth category. Let $M_0$ be an $m$-dimensional closed and connected manifold diffeomorphic to the obtained manifold. We have a special generic map $f_0$ so that $W_{f_0}$ is diffeomorphic to a manifold represented as a boundary connected sum of $|J|$ manifolds in $\{S^{k_j} \times D^{n-k_j}\}_{j \in J}$ just before. Here the boundary connected sum is considered in the smooth category. Furthermore, here we can do in such a way that $N={\mathbb{R}}^n$, that $f_0 {\mid}_{S(f_0)}$ is an embedding and that the two bundles of Proposition \ref{prop:1} (\ref{prop:1.3}) are trivial bundles.

\end{Ex}
\begin{Prop}[E. g. \cite{saeki}]
\label{prop:2}
Let $m>n \geq 1$ be integers. 
A special generic map $f:M \rightarrow N$ on an $m$-dimensional closed and connected manifold $M$ into an $n$-dimensional manifold $N$ with no boundary always enjoys the following properties. We abuse the notation of Proposition \ref{prop:1} such as $q_f:M \rightarrow W_f$.
\begin{enumerate}
	\item We have an {\rm (}$m+1${\rm )}-dimensional compact and connected topological manifold {\rm (PL)} manifold $W$ bounded by $M$, {\rm (}resp. where we discuss in the PL category{\rm )}, and which collapses to an $n$-dimensional compact and smooth manifold $W_f$ in the situation of Proposition \ref{prop:1} and identified with a suitable CW subcomplex {\rm (}resp. subpolyhedron{\rm )} of $W$. If $m-n=1,2,3$ in addition, then $W$ can be chosen as a smooth manifold and $r$ as a smooth map.                                                                                                 
	\item The canonical inclusion $i_M:M \rightarrow W$ and a continuous {\rm (}resp. PL{\rm )} map $r_f:W \rightarrow W_f$ giving a collapsing to $W_f$ satisfy the relation $q_f=r \circ i_M$.
	\item If $M$ is oriented, then $W$ can be chosen as an oriented manifold. If $m-n=1,2,3$ in addition, then $W$ can be chosen as a smooth manifold and $r$ as a smooth map as just before.
    \end{enumerate}
\end{Prop}
We give some remarks on Proposition \ref{prop:2} as Remark \ref{rem:1} in the following.
\begin{Rem}
	\label{rem:1}
Let $m > n \geq 1$ be integers again. Suppose that we have a smooth immersion $\bar{f_0}:W_{f_0} \rightarrow N$ of an $n$-dimensional compact and connected smooth manifold $W_{f_0}$ into an $n$-dimensional manifold $N$ with no boundary. We have a special generic map $f_0:M_0 \rightarrow N$ on a suitable $m$-dimensional closed and connected manifold $M_0$ into $N$ enjoying the three properties in Proposition \ref{prop:1} and we can have an {\rm (}$m+1${\rm )}-dimensional compact and connected smooth manifold $W_0$ and a smooth map $r_0:W_0 \rightarrow W_{f_0}$ as in Proposition \ref{prop:2} where $f:M \rightarrow N$, $W$ and $r:W \rightarrow W_f$ are replaced by $f_0:M_0 \rightarrow N$, $W_0$ and $r_0:W_0 \rightarrow W_{f_0}$ for example. Furthermore, we can do in such a way that the bundles in Proposition \ref{prop:1} (\ref{prop:1.3}) are trivial bundles here.

For more generalized versions of Proposition \ref{prop:2}, see \cite{saekisuzuoka} and see papers \cite{kitazawa0.1,kitazawa0.2,kitazawa0.3} by the author.

Note also that we concentrate on the case $m-n=1$ in main ingredients and we discuss problems in the smooth category in situations of Proposition \ref{prop:2} in main ingredients of our paper.

\end{Rem}

\begin{Thm}
\label{thm:1}
Let $m>0$ be an integer. We have the following conditions for closed manifolds to admit special generic maps into fixed Euclidean spaces where connected sums and boundary connected sums are taken in the smooth category. We also abuse the notation in Propositions \ref{prop:1} and \ref{prop:2} such as the manifold $W_f$.
\begin{enumerate}
\item {\rm (\cite{saeki})}
\label{thm:1.1}
Let $m \geq 2$. An $m$-dimensional closed and connected manifold admits a special generic map into ${\mathbb{R}}^2$ if and only if either of the following cases holds.
\begin{enumerate}
\item $m \neq 4$ and $M$ is homeomorphic to the $m$-dimensional unit sphere.
\item $m=4$ and $M$ is an $m$-dimensional standard sphere.
\item $M$ is diffeomorphic to a manifold represented as a connected sum of the total spaces of smooth bundles over a circle or the $1$-dimensional unit sphere whose fibers are diffeomorphic to homotopy spheres which are not $4$-dimensional exotic spheres.
\end{enumerate}
Furthermore, for a special generic map $f$ here, $W_f$ is diffeomorphic to a boundary connected sum of finitely many copies of $S^1 \times D^1$. 
\item {\rm (\cite{saeki})}
\label{thm:1.2}
Let $m=4,5$. An $m$-dimensional closed and simply-connected manifold admits a special generic map into ${\mathbb{R}}^3$ if and only if either of the following two holds.
\begin{enumerate}
\item $M$ is an $m$-dimensional standard sphere.
\item $M$ is diffeomorphic to a manifold represented as a connected sum of the total spaces of smooth {\rm (}linear{\rm )} bundles over the $2$-dimensional unit sphere whose fibers are $3$-dimensional standard spheres.
\end{enumerate}
Furthermore, for a special generic map $f$ here, $W_f$ is diffeomorphic to a boundary connected sum of finitely many copies of $S^2 \times D^1$. 
\item {\rm (\cite{nishioka})}
\label{thm:1.3}
Let $m=5$. An $m$-dimensional closed and simply-connected manifold admits a special generic map into ${\mathbb{R}}^4$ if and only if either of the following two holds.
\begin{enumerate}
\item $M$ is an $m$-dimehsional standard sphere.
\item $M$ is diffeomorphic to a manifold represented as a connected sum of the total spaces of smooth {\rm (}linear{\rm )} bundles over the $2$-dimensional unit sphere whose fibers are $3$-dimensional standard spheres.
\end{enumerate}
\end{enumerate}
\end{Thm}

For studies on special generic maps, see also \cite{burletderham, furuyaporto} as pioneering ones, for example.

\section{The proof of Main Theorems.}
\label{sec:3}
\subsection{Complex projective spaces whose dimensions are greater than $1$ do not admit special generic maps into the Euclidean spaces whose codimensions are $-1$.}
\label{subsec:3.1}
${\mathbb{C}P}^k$ denotes the $k$-dimensional complex projective space, which is a $k$-dimensional complex space and $2k$-dimensional smooth manifold.
\begin{Thm}
\label{thm:3}
Let $m=2k>2$ be an even integer greater than $2$.
${\mathbb{C}P}^k$ does not admit special generic maps into ${\mathbb{R}}^{m-1}$.
\end{Thm}
\subsection{Fundamental algeberaic topology.}
\label{subsec:3.2}
We introduce notation and terminologies on {\it homology groups}, {\it cohomology groups} and {\it rings} and {\it homotopy groups}. For elementary or advanced theory, consult \cite{hatcher} for example.

Let $A$ be a commutative ring.
Let $(X,X^{\prime})$ be a pair of topological spaces satisfying $X^{\prime} \subset X$ where $X^{\prime}$ can be the empty set. 

The {\it $k$-th homology group} ({\it cohomology group}) whose {\it coefficient ring} is $A$ is denoted by $H_k(X,X^{\prime};A)$ (resp. $H^k(X,X^{\prime};A)$). If $A$ is isomorphic to $\mathbb{Z}$ (resp. $\mathbb{Q}$), then we add "{\it integral}" (resp. "{\it rational}") after "$k$-th". If $X^{\prime}$ is empty, then we may omit ",$X^{\prime}$" in the notation and the $k$-th homology group (cohomology group) of the pair $(X,X^{\prime})$ is also called the {\it $k$-th} {\it homology group} (resp. {\it cohomology group}) of $X$. We can define integral and rational ones similarly.
({\it $j$-th}) {\rm (}{\it co}{\rm )}{\it homology classes} of $(X,X^{\prime})$ (or $X$) whose {\it coefficient ring} is $A$ are elements of the ({\it $j$-th}) (resp. co)homology groups whose coefficient rings are $A$. We can define {\it integral} and {\it rational} ({\it  co}){\it homology classes} similarly.

${\pi}_k(X)$ denotes the {\it $k$-th homotopy group} of a topological space $X$.

Let $(X_1,{X_1}^{\prime})$ and $(X_2,{X_2}^{\prime})$ be pairs of topological spaces satisfying ${X_1}^{\prime} \subset X_1$ and ${X_2}^{\prime} \subset X_2$ where the second topological space of each pair can be empty as before. For a continuous map $c:X_1 \rightarrow X_2$ satisfying $c({X_1}^{\prime}) \subset {X_2}^{\prime}$, $c_{\ast}:H_{\ast}(X_1,{X_1}^{\prime};A) \rightarrow H_{\ast}({X_2},{X_2}^{\prime};A)$ and $c^{\ast}:H^{\ast}({X_2},{X_2}^{\prime};A) \rightarrow H^{\ast}(X_1,{X_1}^{\prime};A)$ denote the homomorphisms defined in a canonical way. For a continuous map $c:X_1 \rightarrow X_2$, $c_{\ast}:{\pi}_k(X_1) \rightarrow {\pi}_k(X_2)$ also denotes the homomorphism between the homotopy groups of degree $k$, which is also defined in a canonical way.

Let $H^{\ast}(X;A)$ denote the direct sum ${\oplus}_j H^j(X;A)$ for every integer $j \geq 0$. The {\it cup product} for a sequence $\{c_j\}_{j=1}^l \subset H^{\ast}(X;A)$ of $l>0$ cohomology classes is important. ${\cup}_{j=1}^l c_j$ denotes this. This gives $H^{\ast}(X;A)$ the structure of a graded commutative algebra over $A$ and this is the {\it cohomology ring} of $X$ whose {\it coefficient ring} is $A$. For a pair $c_1,c_2 \in H^{\ast}(X;A)$, we also use the notation such as $c_1 \sqcup c_2$ for example. We define {\it integral} and {\it rational} {\it cohomology rings} similarly.

The {\it fundamental class} of a compact, connected and oriented manifold $Y$ is the ($\dim Y$)-th homology class, defined as the generator of the group $H_{\dim Y}(Y,\partial Y;A)$, isomorphic to $A$, being compatible with the orientation.

Let $i_{Y,X}:Y \rightarrow X$ be a smooth, PL or topologically flat embedding mapping the boundary $\partial Y$ into the boundary $\partial X$ and the interior ${\rm Int}\ Y$ into the interior ${\rm Int}\ X$. In short, $Y$ is embedded  {\it properly} in $X$.
Let $h \in H_{\dim Y}(X,\partial X;A)$. If the value of the homomorphism ${i_{Y,X}}_{\ast}$ induced by the embedding $i_{Y,X}:Y \rightarrow X$ at the fundamental class of $Y$ is equal to $h$, then $h$ is said to be {\it represented} by the (oriented) submanifold $Y$. 

$A$ may be $\mathbb{Z}/2\mathbb{Z}$, the field of order $2$, or more generally, the commutative ring consisting of elements whose orders are at most $2$. In this case, we may ignore orientations. 

We also need the notions of {\it Poincar\'e duals} in the following four.
\begin{itemize}
	\item The {\it Poincar\'e dual} ${\rm PD}(c_{\rm h}) \in H^{\dim Y-j}(Y,\partial Y;A)$ to a homology class $c_{\rm h} \in H_{j}(Y;A)$.
	\item The {\it Poincar\'e dual} ${\rm PD}(c_{\rm c}) \in H_{\dim Y-j}(Y,\partial Y;A)$
	to a cohomology class $c_{\rm c} \in H^{j}(Y;A)$.
	\item The {\it Poincar\'e dual} ${\rm PD}(c_{\rm h}) \in H^{\dim Y-j}(Y;A)$ to a homology class $c_{\rm h} \in H_{j}(Y,\partial Y;A)$.
	\item The {\it Poincar\'e dual} ${\rm PD}(c_{\rm c}) \in H^{\dim Y-j}(Y;A)$ to a cohomology class $c_{\rm c} \in H^{j}(Y,\partial Y;A)$.
\end{itemize} 
Related to this, Poincar\'e duality theorem for the manifold $Y$ is also important.

We can define the {\it dual} to a homology class of a basis of a free submodule of a homology group of $(X,X^{\prime} \subset X)$ where we consider modules over $A$.
Note that $X^{\prime}$ may be empty. It is a uniquely defined cohomology class in the cohomology group of the pair $(X,X^{\prime})$ in such a way that the degrees of the homology class and the cohomology class are same.
\subsection{Proving Theorem \ref{thm:3} and Main Theorems.}
\label{subsec:3.3}
Theorem \ref{thm:3} with the following theorem yields Main Theorem \ref{mthm:1}.
\begin{Thm}[\cite{kitazawa1}]
\label{thm:4}
		Let $m>n \geq 1$ and $l>0$ be integers.	Let $A$ be commutative ring.
		Let $M$ be an $m$-dimensional closed and connected manifold.
	 We suppose the existence of a sequence $\{a_j\}_{j=1}^l \subset H^{\ast}(M;A)$ of cohomology classes satisfying the following two conditions. 
		\begin{itemize}
			\item The cup product ${\cup}_{j=1}^l a_j$ is not the zero element.
			\item The degree of each cohomology class in the sequence $\{a_j\}_{j=1}^l$ is smaller than or equal to $m-n$ and the sum of the degrees of these $l$ cohomology classes is greater than or equal to $n$.
		\end{itemize}
		Then $M$ does not admit special generic maps into any $n$-dimensional connected manifold which is not compact and which has no boundary.
	\end{Thm}
	We review the original proof in \cite{kitazawa1}. 
	An essential ingredient of the proof has been shown in \cite{kitazawa0} first. This is also reviewed in \cite{kitazawa2,kitazawa3,kitazawa4}.
	
	We omit precise expositions on {\it handles} and their {\it indice} for polyhedra including PL manifolds.
	\begin{proof}[Reviewing the proof of Theorem \ref{thm:4}]
		We abuse the notation in Propositions \ref{prop:1} and \ref{prop:2} for example.
		
		Suppose that $M$ admits a special generic map into an $n$-dimensional connected manifold $N$ which is non-closed and which has no boundary. We can take an ($m+1$)-dimensional compact and connected PL manifold $W$ in Proposition \ref{prop:2}: Remark \ref{rem:1} says that we consider the case $m-n=1,2,3$ only and discuss problems in the smooth category and and we may restrict the case as this. 
		
		$W_f$ is an $n$-dimensional compact and connected manifold smoothly immersed into the connected and non-closed manifold $N$ with no boundary via $\bar{f}:W_f \rightarrow N$. As a result it is (simple) homotopy equivalent to an ($n-1$)-dimensional compact and connected polyhedron. 
		$W$ is (simple) homotopy equivalent to $W_f$. We can see that $W$ is a PL manifold obtained by attaching handle to $M \times \{0\} \subset M \times [-1,0]$ whose indices are greater than $(m+1)-{\dim W_f}=m-n+1$ where the boundary $\partial W=M$ can be identified with $M \times \{-1\} \subset M \times [-1,0]$ naturally.
		We can take a unique cohomology class $b_j \in H^{\ast}(W;A)$ satisfying the relation $a_j={i_M}^{\ast}(b_j)$. $W$ is (simple) homotopy equivalent to an ($n-1$)-dimensional polyhedron. This means that the cup product ${\cup}_{j=1}^l a_j$ is the zero element. This contradicts the assumption. This completes the proof.
	\end{proof}

We show Theorem \ref{thm:3} (Main Theorem \ref{mthm:1}) and Main Theorem \ref{mthm:2}
\begin{proof}[A proof of Theorem \ref{thm:3} {\rm (}Main Theorem \ref{mthm:1}{\rm )} and Main Theorem \ref{mthm:2}]
	We prove results through several STEPs. \\
\ \\
STEP 1 Preparations.

We prove theorems by adopting methods in \cite{kitazawa2,kitazawa3,kitazawa4,kitazawa5} together with fundamental theory in \cite{saeki} for example. However, we adopt shorter expositions or different ones.
We abuse the notation in Propositions \ref{prop:1} and \ref{prop:2} such as a special generic map $f:M \rightarrow N:={\mathbb{R}}^{n}={\mathbb{R}}^{m-1}$, the surjective smooth map $q_f:M \rightarrow W_f$ onto an ($m-1$)-dimensional compact and connected manifold $W_f$ smoothly immersed into ${\mathbb{R}}^{m-1}$ via a smooth immersion $\bar{f}:W_f \rightarrow {\mathbb{R}}^n$, an ($m+1$)-dimensional compact and connected smooth manifold $W$ with a smooth map $r:W \rightarrow W_f$. Assume the existence of such a map $f:M \rightarrow {\mathbb{R}}^n$.
	
The argument on handles to obtain $W$ from $M$ or $M \times [-1,0]$ implies that $M$, $W$ and $W_f$ are simply-connected. The canonical inclusion $i_M:M \rightarrow W$ and $r:W \rightarrow W_f$ induce isomorphisms ${i_{M}}_{\ast}:{\pi}_1(M) \rightarrow {\pi}_1(W)$ and ${q_f}_{\ast}:{\pi}_1(W) \rightarrow {\pi}_1(W_f)$.

For $r:W \rightarrow W_f$, we have the following objects and properties.
\begin{itemize}
	\item There exists an ($m+1$)-dimensional compact and connected smooth manifold $\tilde{W} \supset W$ which is cornered. 
	\item $\tilde{r}:\tilde{W} \rightarrow W_f$ is regarded as the projection of a bundle over $W_f$ whose fiber is (PL) homeomorphic to the unit disk of dimension $m-n+1=m-(m-1)+1=2$ and whose structure group consists of (resp. PL) homeomorphisms of the disk.
	\item $\tilde{W}$ collapses to $W$, which is regarded as a CW subcomplex (resp. subpolyhedron) of $\tilde{W}$, and these two manifolds are (resp. PL) homeomorphic.
	\item ${\tilde{r}}^{-1}(W_f-{\rm Int}\ N(\partial W_f))=r^{-1}(W_f-{\rm Int}\ N(\partial W_f))$ in Proposition \ref{prop:2}.
	\item The preimage of $W_f-{\rm Int}\ N(\partial W_f)$ by $r$ in Proposition \ref{prop:2} is unchanged in the collapsing of $\tilde{W} \supset {\tilde{r}}^{-1}(W_f-{\rm Int}\ N(\partial W_f))=r^{-1}(W_f-{\rm Int}\ N(\partial W_f))$ to $W \supset r^{-1}(W_f-{\rm Int}\ N(\partial W_f))={\tilde{r}}^{-1}(W_f-{\rm Int}\ N(\partial W_f))$. The restrictions of $r$ and $\tilde{r}$ to these subsets agree.
\end{itemize}
STEP 2 Supposing that the singular set $S(f)$ is not connected and a proof of the fact that this is connected in the case $k_1>2$.  \\
Suppose that the singular set $S(f)$ is not connected. Then there exist a $1$-dimensional compact and connected manifold $L$ smoothly and properly embedded into $W_f$ and two  distinct connected components $C_1$ and $C_2$ of the boundary $\partial W_f$ and $L \bigcap C_j$ is not empty for $j=1,2$. In other words, the boundary is embedded in the subset $C_1 \bigcup C_2$ of the boundary $\partial W_f$ and the interior is embedded in the interior ${\rm Int}\ W_f$. The preimage ${q_f}^{-1}(L)$ can be regarded as a $2$-dimensional standard sphere. It is regarded as the $2$-dimensional manifold of the domain of a Morse function for the Reeb's theorem. The ($m-2$)-th integral homology class represented by $S(f)$ is regarded as the Poincar\'e dual to the integral cohomology class regarded as the dual to a homology class represented by the $2$-dimensional sphere if the rank of $H^2(M;\mathbb{Q})$ is not $0$. Note that this is the case of $k_1=2$. More precisely, in this case the rank is $1$. We can define the dual to the 2nd integral homology class in a suitable way since the ranks of $H_2(M;\mathbb{Z})$ and $H_2(M;\mathbb{Q})$ are $1$ here. \\
\ \\
STEP 3 A proof of the fact that in the case $k_1>2$ the rank of the integral homology group $H_2(W_f;\mathbb{Z})$ is $0$.

Suppose that the rank of $H_2(W_f;\mathbb{Z})$ is not $0$. By the structure of the map, there exists a smoothly embedded copy ${S_f}^{2}$ of the $2$-dimensional unit sphere in $M$ enjoying the following properties.
\begin{itemize}
	\item $q_f({S_f}^2) \bigcap \partial W_f$ is a finite set.
	\item The value of the homomorphism ${q_f}_{\ast}$ at the 2nd integral homology class represented by the sphere is some element of $H_2(W_f;\mathbb{Z})$ whose order is infinite. This is seen as a homology class represented by a $2$-dimensional sphere smoothly embedded in ${\rm Int}\ W_f$. ${S_f}^2$ can be seen as a kind of variants of the images of sections of the bundle over the sphere in ${\rm Int}\ W_f$ whose projection is given by the restriction of $q_f$.
	\end{itemize}  
If $k_1>2$, then we can see that the rank of $H_2(M;\mathbb{Z})$ must be not $0$. This is a contradiction. \\
\ \\
STEP 4 A proof of the fact that in the case $k_1=2$ $S(f)$ is connected.

We use contradiction as in STEP 3. Suppose that $S(f)$ is not connected in the case $k_1=2$. The rank of $H_2(W_f;\mathbb{Z})$ must be $0$ due to a related argument in STEP 2 and by the assumption that the ranks of $H_2(M;\mathbb{Q})$ and $H_2(M;\mathbb{Z})$ must be at most $1$. This means that the two bundles in Proposition \ref{prop:1} (\ref{prop:1.3}) must be trivial. This is due to well-known classifications of bundles whose fibers are circles. A smooth bundle whose fiber is a circle is well-known to be regarded as a linear bundle and such bundles over a fixed space which is a CW complex (or a space of some suitable wider class of topological spaces) are classified by the 2nd integral cohomology classes of the base space. See also \cite{milnorstasheff,steenrod} for example. 

We can see that the so-called {\it self-intersection} of $S(f)$ can be taken as the empty set. In other words, we can consider a smooth isotopy moving the original $S(f)$ to a place apart from the original place by the structures of the bundles in Proposition \ref{prop:1} (\ref{prop:1.3}). STEP 2 implies that the square $u_0 \cup u_0$ of $u_0$ is the zero element of $H^4(M;\mathbb{Q})$ due to Poincar\'e duality for $M$. 
We can also see that even if we discuss the case where the coefficient ring is $\mathbb{Z}$, then the cup product is the zero element of $H^4(M;\mathbb{Z})$. The assumption that $S(f)$ is not connected is shown to be false.
\\
\ \\
STEP 5 A proof of the fact that in the case $k_1>2$ the ranks of the integral homology groups $H_j(W_f;\mathbb{Z})$ and $H_{j-1}(W_f,\partial W_f;\mathbb{Z})$ are $0$ for $2 \leq j \leq k_1-1$.

Suppose that the rank of the integral homology group $H_j(W_f;\mathbb{Z})$ is not $0$ here. STEP 3 and the classification of smooth bundles whose fibers are circles, reviewed in STEP 4, imply that the bundles of Proposition \ref{prop:1} (\ref{prop:1.3}) are trivial. This implies that the rank of the integral homology group $H_j(M;\mathbb{Z})$ is not $0$. This is a contradiction. Suppose that the rank of the integral homology group $H_{j-1}(W_f,\partial W_f;\mathbb{Z})$ is not $0$ here. We consider the bundle whose projection is $\tilde{r}:\tilde{W} \rightarrow W_f$, discussed in STEP 1. 
For each element of $H_{j-1}(W_f,\partial W_f;\mathbb{Z})$, we can consider a kind of {\it prism operators} or {\it Thom classes} to obtain a unique element of $H_{(j-1)+(m-n+1)}(W;\mathbb{Z})=H_{m-n+j}(W;\mathbb{Z})=H_{j+1}(W;\mathbb{Z})$ and consider the boundary to obtain a unique element of $H_{j}(M;\mathbb{Z})$ in the context of arguments on the homology groups. Note that $n=m-1$ and that $\tilde{W}$ is (PL) homeomorphic to and collapses to $W$ here for example. This gives a monomorphism and implies that the rank of the integral homology group $H_{j}(M;\mathbb{Z})$ is not $0$. This is also a contradiction. \\
\ \\
STEP 6 A proof of the fact that for any integer $k_1 \geq 2$ the ranks of the integral homology groups $H_{k_1}(W_f;\mathbb{Z})$ and $H_{k_1-1}(W_f,\partial W_f;\mathbb{Z})$ are $1$ and $0$ respectively.

$S(f)$ has been shown to be connected in STEPs 2 and 4. For example Proposition \ref{prop:1} (\ref{prop:1.2}) implies that $\partial W_f$ is connected. We have a homology exact sequence for $W_f$ as 
$$H_{1}(W_f;\mathbb{Z}) \rightarrow H_{1}(W_f,\partial W_f;\mathbb{Z})
\rightarrow H_{0}(\partial W_f;\mathbb{Z}) \rightarrow H_{0}(W_f;\mathbb{Z})
$$
and $H_{1}(W_f,\partial W_f;\mathbb{Z})$ is the trivial group by the facts that $W_f$ is simply-connected and that $\partial W_f$ is connected. This allows us to consider the case $k_1>2$ only in proving that the rank of $H_{k_1-1}(W_f,\partial W_f;\mathbb{Z})$ is $0$.

Note that in STEP 6 (and throughout the present proof), the fact that $M$, $W$ and $W_f$ are simply-connected is important to apply Poincar\'e duality for any coefficient ring of the (co)homology groups.  

Suppose that the rank of $H_{k_1-1}(W_f,\partial W_f;\mathbb{Z})$ is greater than $1$. By the canonical way presented in STEP 5, we can obtain elements of $H_{k_1-1+m-n}(M;\mathbb{Z})=H_{k_1}(M;\mathbb{Z})$. We can also see that the ranks of $H_{k_1}(M;\mathbb{Z})$ and $H_{k_1}(M;\mathbb{Q})$ are greater than $1$. This is a contradiction. Suppose that the rank of $H_{k_1-1}(W_f,\partial W_f;\mathbb{Z})$ is $1$. Similarly, we can see that the ranks of $H_{k_1}(M;\mathbb{Z})$ and $H_{k_1}(M;\mathbb{Q})$ are at least $1$ and that the ranks are $1$ by the assumption. Consider the Poincar\'e dual to the dual to a generator of $H_{k_1-1}(W_f,\partial W_f;\mathbb{Q})$, whose rank is $1$ of course and the resulting element is an element of $H_{m-k_1}(W_f;\mathbb{Q})$ and a generator of this vector space over $\mathbb{Q}$. We may regard that this is an element of a generator of $H_{m-k_1}(W_f;\mathbb{Z})$, whose rank is $1$ of course. The bundles in Proposition \ref{prop:1} (\ref{prop:1.3}) are trivial by the condition $k_1>2$ and STEP 3.
By considering a section of the trivial bundle, we have the Poincar\'e dual to the dual to a generator of $H_{k_1}(M;\mathbb{Q})$. 
 This means that the square of a generator of $H^{k_1}(M;\mathbb{Q})$ or the cup product of two same element here is the zero element. This is a contradiction. The ranks of $H_{k_1-1}(W_f,\partial W_f;\mathbb{Z})$ and $H_{k_1-1}(W_f,\partial W_f;\mathbb{Q})$ are shown to be
$0$.

Suppose that the rank of $H_{k_1}(W_f;\mathbb{Z})$ is greater than $1$. Arguments on sections or variants of bundles over spheres or more general subspaces of $W_f$ in STEP 3 and STEP 5 imply that the rank of $H_{k_1}(M;\mathbb{Z})$ is greater than $1$. This is a contradiction. Suppose that the rank of $H_{k_1}(W_f;\mathbb{Z})$ is $0$. The rank of $H_{k_1}(W_f;\mathbb{Q})$ is $0$ of course then. We have a homology exact sequence for $W$ as
$$H_{k_1+1}(W;\mathbb{Q}) \rightarrow H_{k_1+1}(W,\partial W;\mathbb{Q})
\rightarrow H_{k_1}(M;\mathbb{Q}) \rightarrow H_{k_1}(W;\mathbb{Q})
$$
and by Poincar\'e duality theorem for $W$, $H_{k_1+1}(W,\partial W;\mathbb{Q})$ and $H_{m-k_1}(W;\mathbb{Q})$ are isomorphic as vector spaces over $\mathbb{Q}$. Since $W$ collapses to $W_f$, they are also isomorphic to $H_{m-k_1}(W_f;\mathbb{Q})$ and Poincar\'e duality theorem for $W_f$ implies that they are also isomorphic to $H_{(m-1)-(m-k_1)}(W_f,\partial W_f;\mathbb{Q})=H_{k_1-1}(W,\partial W_f;\mathbb{Q})$. The rank of $H_{k_1-1}(W_f,\partial W_f;\mathbb{Z})$ is shown to be
$0$ for $k_1 \geq 2$ and this implies that the rank of $H_{k_1}(M;\mathbb{Q})$ is $0$ from the related argument in STEP 5. This is a contradiction. We have that the ranks of $H_{k_1}(W_f;\mathbb{Q})$ and $H_{k_1}(W_f;\mathbb{Z})$ are $1$.
 \\
\ \\
STEP 7 Investigating the ($k_2-1$)-th power ${\cup}_{j=1}^{k_2-1} {{u_0}}$ of a generator $u_0 \in H^{k_1}(M;\mathbb{Q})$ and showing that this is the zero element and that this is a contradiction.

The ranks of $H_{k_1}(W_f;\mathbb{Z})$ and $H_{k_1}(W_f;\mathbb{Q})$ have been shown to be $1$. 
We consider a generator of a free subgroup of $H_{k_1}(W_f;\mathbb{Z})$ whose rank is $1$, the dual to this, which we can define in a suitable way, and the Poincar\'e dual to the integral cohomology class, which is an element of the ($m-k_1-1$)-th integral homology group $H_{m-k_1-1}(W_f;\mathbb{Z})$. We canonically have a unique element of $H_{m-k_1}(M;\mathbb{Z})$ by the method in STEP 5. We may regard that this is also the Poincar\'e dual ${\rm PD}(u_0) \in H_{m-k_1}(M;\mathbb{Q})$ to $u_0$.

We consider intersection theory and Poincar\'e duality for $W_f$ and $M$. We can consider (variants of) {\it generic} intersections for copies of the homology class ${\rm PD}(u_0)$. Here we may discuss where the coefficient ring is $\mathbb{Z}$. Remember that $m=k_1k_2$ and that the dimension of $W_f$ is $m-1$.
 We have the relation $$(m-k_1-1)+(k_2-2)\{(m-k_1-1)-(m-1)\}=(m-k_1-1)-k_1(k_2-2)=k_1-1$$
and this relation enables us to regard that the Poincar\'e dual to the ($k_2-1$)-th power ${\cup}_{j=1}^{k_2-1} u_0$ of $u_0$ is represented by the element obtained by the method in STEP 5 from an element of $H_{k_1-1}(W_f,\partial W_f;\mathbb{Q})$, which is the trivial vector space. The ($k_2-1$)-th power is the zero element of $H^{k_1(k_2-1)}(M;\mathbb{Q})$. This is a contradiction. \\
\ \\
This completes our proof.
\end{proof}
\begin{proof}[A proof of Main Theorem \ref{mthm:1}]
Theorem \ref{thm:4} and Main Theorem \ref{mthm:2} completes the proof for the case $n \neq m$, For the case $m=n$, we can know from the theory of Eliashberg \cite{eliashberg} and well-known facts on {\it Pontrjagin classes} of complex projective spaces, explained about in \cite{milnorstasheff} for example. 
	\end{proof}
\begin{Ex}
	The $k_2$-dimensional quaternionic projective space is for $k_1=4$ in Main Theorem \ref{mthm:2}. Its ($j$-th) integral homology group is always free for any integer $j$.
	\cite{fowlersu} presents some meaningful examples for Main Theorem \ref{mthm:2} more. See Theorem C there for example.
\end{Ex}

We close our paper by a remark related to a fact discussed by the author first in \cite{kitazawa5} and another several remarks.
\begin{Rem}
	We consider a situation in Remark \ref{rem:1}. 
	
	Assume that a smooth immersion $\bar{f_0}:W_{f_0} \rightarrow N$ is given similarly. According to \cite{kitazawa5}, we have a special generic map $f_0:M_0 \rightarrow N$ on a suitable $m$-dimensional closed and connected manifold $M_0$ enjoying the following two properties.
	\begin{enumerate}
		\item The two bundles in Proposition \ref{prop:1} (\ref{prop:1.3}) are trivial bundles where the notation is as in Remark \ref{rem:1}.
		\item There exists a special generic map $f_k:M_0 \rightarrow N \times {\mathbb{R}}^k$ such that the composition of $f_k$ with the canonical projection to $N$ is the original map $f_0$ for any integer $1 \leq k \leq m-n$. 
	\end{enumerate} 

More precisely, in Main Theorem 1 of \cite{kitazawa5}, a special generic map $f:M \rightarrow N$ on an $m$-dimensional closed and connected manifold $M$ enjoying some good topological properties is given and $f_0$ is given as a map on the manifold $M_0$, obtained as a manifold whose cohomology ring is isomorphic to that of $M$.

It is an open problem whether $M$ in Main Theorem \ref{mthm:2} admits special generic maps into ${\mathbb{R}}^n$ for $m-k_1<n\leq m$ or $n=m$ where we may regard the case $n=m$ is another problem. Note also that it does not admit ones into ${\mathbb{R}}^n$ for $1 \leq n \leq m-k_1$ by Theorem \ref{thm:4}.
We can also say that the special generic maps are not ones presented here even if partial answers to the open problem are given affirmatively.

We also note that real projective spaces admit no special generic maps in most cases (\cite{kitazawa1}). These cases are also suitable for exercises on Theorem \ref{thm:4}.
\end{Rem}
\begin{Rem}
	In the present paper, we may drop conditions on homology and homotopy groups for example. For example, we concentrate on studies on closed and simply-connected manifolds whose dimensions are at least $5$, which are central objects in classical algebraic topology and differential topology. As pioneers, we and the author are trying to understand these manifolds geometric and constructive ways via special generic maps and more general good smooth maps such as so-called {\it fold} maps, generalizing the class of Morse functions and special generic maps. Related papers are \cite{kitazawa1,kitazawa2,kitazawa3,kitazawa4} for example. 
	For our remark here, see also Remark 4 of \cite{kitazawa2} for example.
	For systematic singularity theory of so-called {\it generic} smooth maps including fold maps, we introduce \cite{golubitskyguillemin} as a related book.
\end{Rem}
\begin{Rem}
	In Main Theorems, it is not meaningful to discuss the case where the Euler number of the manifold is odd. For example, the case where $k$ is even in Main Theorem \ref{mthm:1} and the case where $k_2$ is even in Main Theorem \ref{mthm:2} are such cases. This is due to the existence theory of fold maps in \cite{eliashberg,eliashberg2}
	
	As another remark, in Main Theorem \ref{mthm:2}, $k_1$ must be even to make the problem meaningful.  
\end{Rem}
\section{Acknowledgments and declarations.}
Naoki Kitazawa is a member of the project supported by JSPS KAKENHI Grant Number JP17H06128 "Innovative research of geometric topology and singularities of differentiable mappings"
(Principal investigator: Osamu Saeki). Our present study is supported by the project.

We declare that data essentially supporting the study are all in the present paper.
\end{document}